\documentclass[a4paper,12pt]{article}
\usepackage[T1]{fontenc}
\usepackage[english]{babel}
\usepackage{amsmath}
\usepackage{amsfonts}
\usepackage{lmodern}
\usepackage{moreverb}
\usepackage{amsthm}
\usepackage{hyperref}
\usepackage[utf8]{inputenc}
\usepackage[a4paper,tmargin=1.4cm,bmargin=1.5cm,
rmargin=1.5cm,lmargin=1.5cm]{geometry}

\newtheorem{theorem}{Theorem}[section]
\newtheorem{Lemme}{Lemma}[section]
\newtheorem{Proposition}{Proposition}[section]
\newtheorem{Corollaire}{Corollary}[section]
\theoremstyle{definition}
\newtheorem{Remarque}{Remark}[section]
\newtheorem{definition}{Definition}[section]
\newtheorem{exemple}{Example}[section]

\begin{document}

\begin{center}

\begin{Huge}

  \textbf{On the generalization of the study of a letter power substitution on modulo-recurrent words}
  
  \end{Huge}
   
 \begin{center}

\end{center}
 
 \end{center}
  $\hspace{3cm}$Barro Moussa$^{1,a}$, Bognini K. Ernest$^{2,a}$, Kientéga Boucaré$^{3,a}$\\
   
   $^{1}$mous.barro@yahoo.com, UFR/Sciences exactes et Appliquées, Université Nazi BONI
   
   $^{2}$ernestk.bognini@yahoo.fr, Centre Universitaire de Kaya, Université Joseph Ki-Zerbo
   
   $^{3}$kientega.boucare@gmail.com, Institut Universitaire Professionnalisant, Université de Dédougou
   
   $^{a}$Burkina Faso\\
   
  \begin{center}
 \begin{abstract}
 \noindent
Let us consider an infinite word and $k\geq 1$ an integer. By steps of k, we substitute a letter of this infinite word by the power of an external letter. The new word obtaining by this process is called $k$ to $k$ substitution of a power letter. After the application of this new notion on modulo-reccurent words and in particular on Sturmian words. We establish the complexity function of those words.
  \\[3mm]
 
 {\textbf{Keywords:} Sturmian words, modulo-recurrent words, $k$ to $k$ substitution, factors, complexity function.
\\[2mm]
 {\textbf{2020 Mathematics Subject Clasification:} 68R15, 11B85, 03D15.
}
}
 \end{abstract}
\end{center}

   \section{Introduction}
 The complexity function, which counts the number of distinct factors of given length in
some infinite word, is often used in the characterization of some families of words $\cite{b6}$. For instance,  
Sturmian words are non-periodic infinite words  with minimal complexity $\cite{b10, b11}$.  Over the past thirty years, Sturmian words have been intensively study. Thus, these investigations has led to numerous characterizations 
and various properties $\cite{b27, b5, CKT, b23, b12}$ on these words. 
 
The notion of $k$ to $k$ insertion of a letter on infinite words was 
introduced in $\cite{b23}$, and widely studied in $\cite{BKT, b24}$. Later, in $\cite{BBT}$, the first two authors introduced the concepts of $k$ to $k$ substitution and $k$ to $k$ exchange of an internal letter on infinite words.

Here, we present a generalization of $k$ to $k$ substitution introduced and studied in $\cite{BBT}$. More precisely, we now need to substitute in steps  $k$ a letter power external to the alphabet used to construct some infinite word namely modulo-recurrent words and in particular Sturmian words. Thus, the new word obtaining by this application is so called \textit{word by $k$ to $k$ substitution} of letter.

This article focuses on the study of the combinatorial properties of words obtained by this notion on modulo-recurrent words in general and for Sturmian words in particular. The paper is organized as follow. In Section 2, we give useful definitions and notations in combinatorics on words. We also recall some properties of Sturmian words and modulo-recurrent words . In Section 3, we introduce the notion of $k$ to $k$ substitution of a letter power in infinite words in Subsection 3.1. Then, in Subsection 3.2, 
We show that there is a relationship between the complexity function of the word modulo-recurrent words and these $k$ to $k$ substitute words. Finally, we establish a same results for the case of Sturmian words in subsection 3.3.


\section{Background}

\subsection{Definitions and notations}

  An alphabet $\mathcal{A}$, is a non empty finite set whose the elements are called letters. A word is  a finite or infinite sequence of  elements over $\mathcal{A}$. The set of finite words over $\mathcal{A}$ is denoted $\mathcal{A}^\ast$
and $\varepsilon$, the empty word. For any $u\in \mathcal{A}^*$, the number 
of letters of $u$ is called length of $u$ and it is denoted $|u|$. Moreover, for any letter $x$ of $\mathcal{A}$, $|u|_x$ is the number of occurrences of $x$ in $u$. A word $u$ of length $n$ written with a unique letter $x$ is simply denoted $u=x^n$.

 Let $u=x_1x_2 \dotsb x_n$ be a word such that $x_i\in \mathcal{A}$, for all $i\in\left\{1,2,\dots , n\right\}$. 
The image of $u$ by the reversal map is the word denoted $\overline{u}$ and defined by   $\overline{u}=x_n\dotsb$ $x_2x_1$.  The word $\overline{u}$ is simply called reversal image of $u$.  A finite word $u$ is called palindrome if $\overline{u}=u$. If $u$ and $v$ are two finite words over $\mathcal{A}$, we have $\overline{uv}$=$\overline{v}   \hspace{0.05cm} \overline{u}$.

 The set of infinite words over $\mathcal{A}$ is denoted $\mathcal{A}^{\omega}$ and we write $ \mathcal{A}^\infty=\mathcal{A}^*\cup \mathcal{A}^{\omega}$, the set of finite and infinite words.  
An infinite word $\textbf{u}$ is said to be eventually periodic if there exist two words $v\in \mathcal{A}^*$ and $w\in \mathcal{A}^+$ 
such that $u=vw^\omega$. If $v =\varepsilon $, then  $u$ is periodic.
 
Let $\textbf{u}\in \mathcal{A}^\infty$  and $w\in \mathcal{A}^*$. The word $w$ is a factor of $u$ if there exist $u_1\in \mathcal{A}^*$ and $\textbf{u}_2\in \mathcal{A}^\infty$ such that $\textbf{u}=u_1w\textbf{u}_2$. The factor $w$ is said to be a prefix (respectively, a suffix) if $u_1$ (respectively, $\textbf{u}_2$) is the empty word.

  A word $\textbf{u}$ is said to be recurrent if each of its factors appears infinitly in $\textbf{u}$.
 A word $\textbf{u}$ is said to be uniformly recurrent if for all integers $n$, there exists an integer $N$ such that any factor of length $N$ in $\textbf{u}$ contains all the factors of length $n$.

A non-empty factor $w$ of $\textbf{u}$, is said to be right (respectively, left) extendable by a letter $x$ in $\textbf{u}$ if 
$wx$ (respectively, $xw$) appears in $\textbf{u}$. The number of right (respectively, left) extensions of $w$, is denoted $\partial^+w$ (respectively, $\partial^-w$). The factor $w$ is said to be right (respectively, left) special in $\textbf{u}$ if $\partial^+w>1$ (respectively, $\partial^-w>1$). If $\partial^+w=2$ (respectively, $\partial^+w=k$), we say that $w$ have two-right (respectively, $k$-right) extensions in $\textbf{u}$.  Idem for left extensions. The factor $w$ is said to be bispecial in $\textbf{u}$ if $w$ is both right and left special.
    
Let $\textbf{u}$ be an infinite word over $\mathcal{A}$. The set of factors of length $n$  in $\textbf{u}$, is written $\text{L}_n(\textbf{u})$ and the set of all factors in $\textbf{u}$ is denoted by $\text{L}(\textbf{u})$. Let $\textbf{u}=x_0x_1x_2 \cdots $, where $x_i\in \mathcal{A}$, $i\geq 0$ be an infinite word and $w$ his factor. Then, $w$ appears in $\textbf{u}$ at the position $l$ if $w=x_{l}x_{l+1} \cdots x_{l+|w|-1}$.

The complexity function of a given infinite word $\textbf{u}$ is the map of $\mathbb{N}$ to $\mathbb{N}^\ast$ defined by $\textbf{P}_\textbf{u}(n) = \#\text{L}_n(\textbf{u})$, where $\#\text{L}_n(\textbf{u})$ designates the cardinal of $\text{L}_n(\textbf{u})$.

 This function is related to the special factors by the relation (see \cite{b8}):
    $$\hspace{0 cm} \textbf{P}_\textbf{u}(n+1)-\textbf{P}_\textbf{u}(n)= \displaystyle\sum_{w\in \text{L}_n(\textbf{u})} (\partial^+(w)-1)=\displaystyle\sum_{w\in \text{L}_n(\textbf{u})} (\partial^-(w)-1). $$
 
We call first difference of the complexity function $\textbf{P}_{\textbf{u}}$, the function defined for any integer $n$ by:
 $$ \textbf{S}_\textbf{u}(n)=\textbf{P}_\textbf{u}(n+1)-\textbf{P}_\textbf{u}(n), \forall n \in \mathbb{N}. $$
  


 Let $\textbf{u}=x_0x_1x_2x_3\cdots $ be an infinite word. The window complexity function of $\textbf{u}$ is the map,  $\textbf{P}^f_{\textbf{u}}$ of $\mathbb{N}$ into $\mathbb{N}^*$, defined by: $$\textbf{P}^f_{\textbf{u}}(n)=\#\left\{x_{kn}x_{kn+1}\cdots x_{n(k+1)-1} : k\geq 0\right\}.$$ 

A morphism $f$ is a map of $\mathcal{A}^*$ into itself such that $f(uv)=f(u)f(v)$, for 
any $u$, $v\in \mathcal{A}^*$. 
  
\subsection{Sturmian words and modulo-recurrent words}
 In this subsection, we recall definitions and some useful properties on Sturmian words and modulo-recurrent words.
\begin{definition}\label{def-stur} An infinite word $\textbf{u}$ over $\mathcal{A}_{2}=\{a,b\}$ is said to be Sturmian if for any integer $n$, $\textbf{P}_\textbf{u}(n)=n+1$.
  \end{definition} 
The well-known Sturmian word in the literature is the famous Fibonacci word. It is generated by the 
morphism  $ \varphi $ defined over $\mathcal{A}_{2}=\{a,b\}$ by:
$$\varphi(a)=ab \ \text{and} \ \varphi(b)=a.$$ 

\begin{definition}
An infinite word $\textbf{u}=x_0x_1x_2 \dotsb$ is said to be modulo-recurrent if, any factor $w$ of $\textbf{u}$ appears in $\textbf{u}$ at all positions modulo $i$, for all $i\geq 1$.

\end{definition}

 \begin{Proposition}$\cite{CKT}$\label{prop2-stur-mod}
Let $\textbf{\emph{u}}\in \mathcal{A}^{\infty}$ such that $\textbf{P}_\textbf{\emph{u}}(n)=(\#\mathcal{A})^n$, for all $n\in \mathbb{N}$. Then, $\textbf{\emph{u}}$ is a modulo-recurrent word.
\end{Proposition} 

 \begin{definition} Let $\textbf{u}$ be a modulo-recurrent word. Then, $\textbf{u}$ is called non-trivial if there exists an integer $ n_0$ such that for all $n\geq n_0$:  $$ \textbf{P}_\textbf{u}(n)<(\#\mathcal{A})^n.$$ 
\end{definition}

\begin{definition}  Let $\textbf{u}$ be an infinite word. Then, $\textbf{u}$ is said to be uniformly modulo-recurrent if it is uniformly recurrent and modulo-recurrent.
\end{definition} 
\begin{definition}
A factor $w$ of some infinite word $\textbf{u}$ is said to be 
 a window factor when it appears in $\textbf{u}$ at a mutiple position its length.
\end{definition}
Let us denote $\text{L}^f_n(\textbf{u})$, the set of $n$-window factors of $u$ for a given factor of length $n$. Thus, his cardinal is $\textbf{P}^f_\textbf{u}(n)$. 

\begin{theorem}\label{stur-mod}$\cite{b23}$
Every Sturmian word is modulo-recurrent.
 
 \end{theorem}
 (i) It is clear that the Sturmian words are non-trivial and uniformly modulo-recurrent.
  
  (ii) The champernowne word is modulo-recurrent but does not uniformly recurrent.\\
  

The modulo-recurrent words can be characterized by their  window complexity as follow:

\begin{theorem}\label{theo-stur-mod}\cite{CKT} A recurrent word $\textbf{\emph{u}}$ is modulo-recurrent 
if and only if $ \textbf{P}^f_\textbf{\emph{u}}(n)= \textbf{P}_\textbf{\emph{u}}(n)$, for all $ n\geq 1$.

\end{theorem}     
\section{Studies on the complexity function of infinite $k$ to $k$ substituted words}
In this section, we introduce firstly  the notion of $k$ to $k$ substitution on infinite words and we show a fundamental result, Lemme \ref{(k+i)q} for further work (Section 3.1). Next, we apply this notion on modulo-recurrent words and we find its complexity function, Theorem \ref{prop-mod} in Section 3.2. Finally, in Section 3.3, we apply this notion to Sturmian words, which are particular modulo-recurrent words and give their complexity function Theorem \ref{cci}.   
\subsection{Notion of a letter power substitutiton on infinite words}
This subsection is devoted to introducing the notion of $k$ to $k$ substitution of a letter power in infinite words. Now, let us consider an infinite word $\textbf{u}$ over $\mathcal{A}$ in the form:
\begin{equation}\label{equa1}
\textbf{u}=x_0w_1x_1w_2x_2w_3x_3 \dotsb w_ix_i \dotsb,
\end{equation}
with $w_i \in \text{L}_k(\textbf{u})$ and $x_i \in \mathcal{A}$, for all $i\in\mathbb{N}$.\\
Let $x$ be an external letter over $\mathcal{A}$ and $l\geq 1$ an integer. Then we have a useful definition for the remainder of the paper.
 \begin{definition}\label{def1-ins-c} We call, word by $k$ to $k$ substitution of a letter power, the word $\textbf{v}=\mathcal{S}_{k}^{x^l}(\textbf{u})$ obtained by substituting the letter $x_{i}$ by $x^{l}$ in $\textbf{u}$. We have:
 \begin{equation}
 \textbf{v}=\mathcal{S}_{k}^{x^l}(\textbf{u})=x^lw_1x^lw_2x^lw_3x^l\dotsb x^lw_ix^l \dotsb.
\end{equation}    
 
                \end{definition}
       
\begin{exemple}\label{expl2-ins}
       
      Let us consider the Fibonacci word $\textbf{f}$ over $\{a,b\}$ whose first letters are given by: $$  \textbf{f}=abaababaabaababaababaababaabaababaa\dotsb. $$ 
      Then, in steps of $k=2$ and for $c\notin \{a,b\}$, we obtain the following word:
       $$S_{2}^{c^3}(\textbf{f})=c^3bac^3bac^3aac^3aac^3abc^3abc^3bac^3bac^3aac^3aac^3abc^3a\dotsb $$   called word by $2$ to $2$ substitution of $c^3$ in $\textbf{f}$.  
\end{exemple}   

\begin{Lemme}\label{(k+i)q} Let $\textbf{\emph{u}}$ be an infinite word over $\mathcal{A}$ and $\textbf{\emph{v}}=S_{k}^{x^l}(\textbf{\emph{u}})$. Then, any factor of $\textbf{\emph{v}}$ of length $(k+l)q$ comes from a factor of $\textbf{\emph{u}}$ of length $(k+1)q$.
\end{Lemme}
\textbf{Proof.} Let us firstly observe that any factor of $\textbf{v}$ of length $(k+l)q$ comes from a factor of $\textbf{u}$ of lengths $(k+1)q-1$, $(k+1)q$ or $(k+1)q+1$. Now, we show that the factors of $\textbf{v}$ of length $(k+l)q$ generated by those of $\textbf{u}$ of lengths $(k+1)q-1$ and $(k+1)q+1$ are also generated by those of $\textbf{u}$ of length $(k+1)q$. Let $v_1$ be a factor of $v$ of length $(k+l)q$ from a factor $u_1$ of $\textbf{u}$. Note that each factor of $\textbf{v}$ of length $k+l$ contains $l$ substituted letters. So, it comes from a factor of $\textbf{u}$ of length $k+1$. Since $|v_1|=(k+l)q$, then $v_1$ contains exactly $lq$ substituted letters. In addition, $|v_1|=(k+l)q=kq+lq$. Hence, $|u_1|=(k+1)q$. \hfill$\square$

 \begin{Remarque}\label{propo-rec} Let $\textbf{u}$ be a modulo-recurrent word over $\mathcal{A}$ and $\textbf{v}=S_{k}^{x^l}(\textbf{u})$. Then, the followings holds:
 \begin{enumerate}
\item the word $\textbf{v}$ is recurrent.
\item the factors of $\textbf{v}$ of length $(k+l)q$ starting and ending with $c$ come also from factors of $\textbf{u}$ of length $(k+1)q-1$ and $(k+1)q+1$.
 \end{enumerate}
\end{Remarque}

\begin{Remarque}\label{prol}
Let $\textbf{v}=S_{k}^{x^l}(\textbf{u})$ for a given infinite word $\textbf{u}$.  Let $u_1$ and $u_2$ be two distinct factors of $\textbf{u}$ neither coming from the same factor nor extendable at left (resp. right) by a letter. By substituting the first (resp. last) letter of $u_1$ and $u_2$ with $x^l$, then we obtain the same factor of $\textbf{u}$.
\end{Remarque}

In the following, $\textbf{u}$ denotes an infinite word on $\mathcal{A}$ and $\textbf{v}=S_{k}^{c^l}(\textbf{u})$ with $c\notin\mathcal{A}$ unless otherwise specified.

 \subsection{ $k$ to $k$ substitution of a letter power in modulo-recurrent words}
Here, $\textbf{u}$ is a modulo-recurrent word. The aim of this subsection is to determine the complexity function of the word $\textbf{v}=S_{k}^{c^l}(\textbf{u})$. We show that there is a relationship between the complexity function of the word $\textbf{v}$ and that of $\textbf{u}$.
\begin{theorem}\label{prop-mod} The complexity function of $\textbf{\emph{v}}$ denoted $\textbf{P}_\textbf{\emph{v}}$ is given by the following relations, for all integer $n\geq 1$ : 
 \begin{enumerate}
\item[(i)] For $n<k+l$, we put $m_0= \min  \left\{l, k\right\}$ and $M_0= \max  \left\{l, k\right\}$. Then, we have:
 
 $$\hspace{0 cm} \textbf{P}_\textbf{\emph{v}}(n) = \left \{
 \begin{array}{l}
 \hspace{0cm}1+2\displaystyle\sum_{i=1}^{n-1} \textbf{P}_\textbf{\emph{u}}(i)+\textbf{P}_\textbf{\emph{u}}(n) \ \text{if} \ n\leq m_0 \hspace{0.95cm}\\
 \hspace{0cm}1+2\displaystyle\sum_{i=1}^{k-1} \textbf{P}_\textbf{\emph{u}}(i)+(n-k+1)\textbf{P}_\textbf{\emph{u}}(k) \ \text{if} \  m_0 < n \leq M_0 \ \text{and}\ k<l\\
 \hspace{-0.2cm}(n-l+1)\textbf{P}_\textbf{\emph{u}}(n-l+1)-2\textbf{S}_\textbf{\emph{u}}(n-l)+2\displaystyle\sum_{i=n-l+1}^{n-1} \textbf{P}_\textbf{\emph{u}}(i)+\textbf{P}_\textbf{\emph{u}}(n) \ \text{if} \ m_0 < n \leq M_0 \ \text{and} \ l<k\\
 (n-l+1)\textbf{P}_\textbf{\emph{u}}(n-l+1)-2\textbf{S}_\textbf{\emph{u}}(n-l)+2\displaystyle\sum_{i=n-l+1}^{k-1} \textbf{P}_\textbf{\emph{u}}(i)+(n-k+1)\textbf{P}_\textbf{\emph{u}}(k) \ \text{if} \ n>M_0.
  
 \end{array}
 \right. $$
 
\item[(ii)] For $n\geq k+l$, let us put respectively  $n=(k+l)q+\alpha$ where $0\leq \alpha <k+l$ for $q\geq 1$ and $\beta=\min  \left\{\alpha, k\right\}$. Then, we have:
 
  $$\hspace{0 cm} \textbf{P}_\textbf{\emph{v}}(n) = \left \{
 \begin{array}{l}
 (k-1)\textbf{P}_\textbf{\emph{u}}((k+1)q)+(l+1)\textbf{P}_\textbf{\emph{u}}((k+1)q-1)\ \text{if} \ \alpha=0 \\
 
 (l-\alpha+1)\textbf{P}_\textbf{\emph{u}}((k+1)q-1)+2\displaystyle\sum_{i=1}^{\beta-1} \textbf{P}_\textbf{\emph{u}}((k+1)q+i)+(|k-\alpha|+1)\textbf{P}_\textbf{\emph{u}}((k+1)q+\beta) \\ \text{if} \ 0<\alpha \leq l\\
 (\alpha-l+1)\textbf{P}_\textbf{\emph{u}}((k+1)q+\alpha-l+1)-2\textbf{S}_\textbf{\emph{u}}((k+1)q+\alpha-l)+2\displaystyle\sum_{i=\alpha-l+1}^{\beta-1} \textbf{P}_u((k+1)q+i)+\\(|k-\alpha|+1)\textbf{P}_\textbf{\emph{u}}((k+1)q+\beta) \ \text{if} \ l<\alpha< k+l. 
 \end{array}
 \right. $$
 \end{enumerate} 
 \end{theorem}

\textbf{Proof.} The demonstration of this result is based respectively on the value $l$ of the power of the external letter and the value $k$ of the substitution step. Thus, we distinguish the following cases.\\
 \noindent\textbf{Case 1.} Assume $k\leq l$. Then, we have $m_{0}=k$ and $M_{0}=l$.
 
\begin{itemize}
\item[(i)] Let us show now the following subcases of this item: \\
\textbf{Case 1.1.} if $n\leq k$ then, any factor of $\textbf{u}$ of length $n$ is also a factor of $\textbf{v}$.
 In addition, some factors of $\textbf{v}$ of length $n$ contain letter power of $c$. So, they come from factors of $\textbf{u}$ of length $n-1$ at most. Then, these factors begin or end with a power of letter $c$ because $n \leq l$. Hence, we obtain:
$$ L_n(\textbf{v})=L_n(\textbf{u}) \cup \{c^n, \ c^{n-j}t, \ tc^{n-j} : \ t\in L_j(\textbf{u}), \ j=1, \ldots, n-1 \}.$$
 
 As a result, we have:
  $$\textbf{P}_\textbf{v}(n)=\textbf{P}_\textbf{u}(n)+\textbf{P}_\textbf{u}(0)+2[\textbf{P}_\textbf{u}(1)+\textbf{P}_\textbf{u}(2)+\dotsb +\textbf{P}_\textbf{u}(n-1)].$$
\textbf{Case 1.2.} if $k<n \leq l$ then, the factors of $\textbf{v}$ of length $n$ contain at least one occurrence of $c$. Moreover, these factors can contain up to $n$ occurrences of $c$ because $n\leq l$. Thus, they come from factors of $\textbf{u}$ of length $k$ at most. Thus, the set of factors of $\textbf{v}$ of length $n$ is given by: 
 
 $$L_n(\textbf{v})=\{c^n,\ c^{n-j}t, \ tc^{n-j} : \ t\in L_j(\textbf{u}),\ j=1, \dotsb , k \} \cup \{ c^{j}wc^{n-k-j} : \ w \in L_k(\textbf{u}), \ j=1, \ldots , n-k-1\}.$$
   
    As a result, we have :
  \begin{align*}
\textbf{P}_\textbf{v}(n)&=\textbf{P}_\textbf{u}(0)+2[\textbf{P}_\textbf{u}(1)+\textbf{P}_\textbf{u}(2)+\dotsb +\textbf{P}_\textbf{u}(k)]+(n-k-1)\textbf{P}_\textbf{u}(k)\\
&= 1+2[\textbf{P}_\textbf{u}(1)+\textbf{P}_\textbf{u}(2)+\dotsb +\textbf{P}_\textbf{u}(k-1)]+(n-k+1)\textbf{P}_\textbf{u}(k).
  \end{align*}
   
 \textbf{Case 1.3.}  if $l<n< k+l$ then, the factors of $\textbf{v}$ of length $n$ contain $l-1$, $l-2$, $\dotsb $, $n-k+1$ and $n-k$ occurrences of $c$. Therefore, these factors come respectively from the factors of $\textbf{u}$ of length $n-l+1$, $n-l+2$, $\dotsb$, $k-1$ and $k$. Thus, we obtain:$\vspace{0.2cm}$
 
 $L_n(\textbf{v})= \left\{c^lx_2\dotsb x_{n-l+1} , \hspace{0.1cm} x_1c^{l}x_3\dotsb x_{n-l+1} ,\dotsb , \hspace{0.1cm} x_1\dotsb c^{l}x_{n-l+1}, \hspace{0.1cm} x_1\dotsb x_{n-l}c^{l} : \hspace{0.1cm} \hspace{0.1cm} x_1\dotsb x_{n-l+1}\in L(\textbf{u}) \right\} \vspace{0.2cm}$
 
$\hspace{1.2cm} \cup \left\{c^{j}t \hspace{0cm},\hspace{0.1cm} tc^{j} :\hspace{0.1cm} t \in L_{n-j}(\textbf{u});\hspace{0.1cm} j=n-k,\ldots, l-1 \right\} \cup \left\{ c^{j}wc^{n-k-j}:  \hspace{0.1cm} w \in L_k(\textbf{u}) ;\hspace{0.1cm} j=1,\ldots , n-k-1  \right\}$.$\vspace{0.2cm}$

 Consequently, we have : 
  \begin{align*}
\textbf{P}_\textbf{v}(n)&=(n-l+1)\textbf{P}_\textbf{u}(n-l+1)-2\textbf{S}_\textbf{u}(n)+2[\textbf{P}_\textbf{u}(n-l+1)+\dotsb +\textbf{P}_\textbf{u}(k)]+(n-k-1)\textbf{P}_\textbf{u}(k)\\
&=(n-l+1)\textbf{P}_\textbf{u}(n-l+1)-2\textbf{S}_\textbf{u}(n)+2[\textbf{P}_\textbf{u}(n-l+1)+\dotsb +\textbf{P}_\textbf{u}(k-1)]+(n-k+1)\textbf{P}_\textbf{u}(k)\\
&=(n-l+1)\textbf{P}_\textbf{u}(n-l+1)-2\textbf{S}_\textbf{u}(n)+(n-k+1)\textbf{P}_\textbf{u}(k)+2\displaystyle\sum_{i=n-l+1}^{k-1} \textbf{P}_\textbf{u}(i).
  \end{align*}
    \item[(ii)] $n\geq k+l$. Then, there are two integers $\alpha$ and $q$ such that $n=(k+l)q+\alpha$ with $q\geq 1$ and $0\leq \alpha < k+l$. Then,  let's reason according to the values of $\alpha$:\\
  
\textbf{Case 1.4. } if $\alpha=0 $ then by Lemma \ref{(k+i)q}, the factors of $\textbf{v}$ of length $n$ come from those of $\textbf{u}$ of length $(k+1)q$. Thus, these factors contain $(l-1)q$ occurrences of the letter $c$. These factors contain $q$ or $q+1$ substitutions. We have the set: 
 $L_n(\textbf{v})=\{w_0c^lw_1c^l\cdots c^lw_q; |w_0|+|w_q|=k, 1\leq |w_0|,|w_q|\leq k-1; w_0x_1w_1x_2\cdots x_qw_q\in L_{(k+1)q}(\textbf{u})\}.$
 
$\hspace{0.6cm}\cup \{ c^jw_1c^l\cdots c^lw_qc^{l-j}; \ j= 0,1,2,\ldots,l-1,l:\ w_1x_2\cdots x_qw_q\in L_{(k+1)q-1}(\textbf{u}) \}$.
 
  Hence, we obtain the following complexity: 
 $$\textbf{P}_\textbf{v}(n)=(k-1)\textbf{P}_\textbf{u}((k+1)q)+(l+1)\textbf{P}_\textbf{u}((k+1)q-1).$$
 \textbf{Case 1.5.} if  $ 1\leq \alpha \leq k$, let $v_1$ be a factor of $\textbf{v}$ of length $n$, generated by a factor $u_1$ of $\textbf{u}$. According to Lemma \ref{(k+i)q}, we obtain $|u_1|=(k+1)q+j$ with $j=0,1, \ldots ,\alpha$.
 
\hspace*{0.5cm} $\bullet $ If $|u_1|=(k+1)q$, then we have $|v_1|_c=(l-1)q+\alpha$. Hence, $v_1$ contains $q-1$ internal substitutions and one substitution at the beginning and one at the end of $v_1$. Thus, $v_1$ can be written as $v_1=c^iw_1c^lw_2\dotsb c^lw_qc^{\alpha+l-i} $ where $i=\alpha,\ldots, l$ with $u_1=x_1w_1x_2w_2\dotsb x_qw_q $ or $u_2=w_1x_1w_2x_2\dotsb w_q x_{q}$ and $|w_i|=k$. Note that the factors $v_1$ is generated by $w_1x_2w_2\dotsb x_qw_q=x_{1}^{-1}u_1=u_2x_{q}^{-1}$ which is of length $(k+1)q-1$. Consequently, $u_1$ produces ($l-\alpha+1$) distinct factors of $v$.

\hspace*{0.5cm} $\bullet$ If $|u_1|=(k+1)q+j$ with $j=1, \dotsb, \alpha-1$, then we have $|v_1|_c=(l-1)q+\alpha-j$. Thus, the factor $v_1$ contains $q$ internal substitutions of $c^l$. Consequently, let us observe that the first substitution of $c^l$ in $u_1$ occurs either at position $j$ or at position $k$. Thus, $u_1$ produces two factors of the form $v_1=w_0c^lw_1c^lw_2\dotsb c^lw_qc^{\alpha-j} $ and $v_2=c^{\alpha-j} t_1c^lt_2\dotsb c^lt_qc^{l}t_{q+1} $ with $u_1=w_0x_1w_1x_2w_2\dotsb x_qw_q =t_1y_1t_2y_2\dotsb t_qy_qt{q+1}$ and $|w_i|=|t_i|=k$ for $i=1, \ldots, q$ and $|w_0|=|t_{q+1}|= j$.

\hspace*{0.5cm} $\bullet$ If $|u_1|=(k+1)q+\alpha$, then we have $|v_1|_c=(l-1)q$. Thus, $v_1$ contains $q$ internal substitutions of $c^l$. Consequently, the first substitution of $c^l$ in $u_1$ is between positions $\alpha$ and $k$. Thus, $v_1$ is of the form $v_1=w_0c^lw_1c^lw_2\dotsb c^lw_q$ with $u_1=w_0x_1w_1x_2w_2\dotsb x_qw_q$ and $|w_0|+|w_q|=k+\alpha$, $\alpha \leq |w_0|$, $|w_q|\leq k$. Hence, $u_1$ produces ($k-\alpha+1$) distinct factors of $\textbf{v}$ of length $(k+l)q+\alpha$. Hence, we obtain:  
\begin{align*}
\textbf{P}_\textbf{v}(n)&=(l-\alpha+1)(\textbf{P}_\textbf{u}((k+1)q)-2\textbf{S}_\textbf{u}(n))+2[\textbf{P}_\textbf{u}((k+1)q+1)+\dotsb +\textbf{P}_\textbf{u}((k+1)q+\alpha-1)]+\\
&(k-\alpha+1)\textbf{P}_\textbf{u}((k+1)q+\alpha)\\
&= (l-\alpha+1)(\textbf{P}_\textbf{u}((k+1)q)-2\textbf{S}_\textbf{u}(n))+2\displaystyle\sum_{i=1}^{\alpha-1} \textbf{P}_\textbf{u}((k+1)q+i)+(k-\alpha+1)\textbf{P}_\textbf{u}((k+1)q+\alpha.
\end{align*}
   
\hspace*{0.5cm} \item For  $k<\alpha\leq l$. Let $v_1$ be a factor of $\textbf{v}$ of length $n$, generated by a factor $u_1$ of $\textbf{u}$. By Lemma \ref{(k+i)q}, we have $|u_1|=(k+1)q+j$ with $j=0,\ldots, k$.
 
\hspace*{1cm} $\bullet $ If $|u_1|=(k+1)q$, then we have $|v_1|_c=(l-1)q+\alpha+1$. Thus, $v_1$ contains $q-1$ internal substitutions and one substitution at the beginning and one at the end of $v_1$. Threfore, $v_1$ is written as $v_1=c^iw_1c^lw_2\dotsb c^lw_qc^{\alpha+l-i} $ where $i=\alpha,\ldots, l$ with $u_1=x_1w_1x_2w_2\dotsb x_qw_q $ or $u_2=w_1x_1w_2x_2\dotsb x_qw_q x_{q}$ and $|w_i|=k$. Note that the factors $v_1$ is generated by the factor $x_1^{-1}u_1=u_2x_q^{-1}=w_1x_2w_2\dotsb x_qw_q $, which is of length $(k+1)q-1$. Thus, each factor of length $(k+1)q-1$ produces $(l-\alpha+1)$ distinct factors of $\textbf{v}$ of length $n$.

\hspace*{1cm} $\bullet$ If $|u_1|=(k+1)q+j$ with $j=1,\ldots, k-1$, then we have $|v_1|_c=(l-1)q+\alpha-j$. Thus, $v_1$ contains $q$ substitutions of $c^l$. Consequently, the first substitution of $c^l$ into $u_1$ occurs either at position $j+1$ or at position $k+1$. Thus, $u_1$ produces two factors of the form $v_1=w_0c^lw_1c^lw_2\dotsb c^lw_qc^{\alpha-j} $ and $v_2=c^{\alpha-j} t_1c^lt_2\dotsb c^lt_qc^{l}t_{q+1} $ with $u_1=w_0x_1w_1x_2w_2\dotsb x_qw_q =t_1y_1t_2y_2\dotsb t_qy_qt{q+1}$ and $|w_i|=|t_i|=k$, $|w_0|=|t_{q+1}|= j$.

\hspace*{1cm} $\bullet$ If $|u_1|=(k+1)q+k$, then we have $|v_1|_c=(l-1)q+\alpha-k$. Thus, $v_1$ contains $q$ substitutions of $c^l$. As a result, $v_1$ is of the form  $v_1=c^{h}w_0c^lw_1c^lw_2\dotsb c^lw_qc^{\alpha-k-h}$ with $h=0,1,\ldots, \alpha-k$ and $u_1=w_0x_1w_1x_2w_2\dotsb x_qw_q$ where $|w_i|=k$; $i=0,1,\ldots, q$.

Consequently, the factor $u_1$ produces ($\alpha-k+1$) distinct factors of $\textbf{v}$ of length $(k+l)q+\alpha$. Hence, we have: 
 \begin{align*}
 \textbf{P}_\textbf{v}(n)&=(l-\alpha+1)\textbf{P}_\textbf{u}((k+1)q-1)+2[\textbf{P}_\textbf{u}((k+1)q+1)+\dotsb +\textbf{P}_\textbf{u}((k+1)q+k-1)]+\\
 &(\alpha-k+1)\textbf{P}_\textbf{u}((k+1)q+k)\\
&=(l-\alpha+1)\textbf{P}_\textbf{u}((k+1)q-1)+2\displaystyle\sum_{i=1}^{k-1} \textbf{P}_\textbf{u}((k+1)q+i)+(\alpha-k+1)\textbf{P}_\textbf{u}((k+1)q+k).
 \end{align*}
 
\textbf{Case 1.6.} if $\alpha>l$. By Lemma \ref{(k+i)q}, the factors of $\textbf{v}$ of length $n$ come from the factors of $\textbf{u}$ of length $(k+1)q+r+1$, $(k+1)q+r+2$, $ \dotsb$, $(k+1)q+k-1$ and $(k+1)q+k$ where $r=\alpha-l$. Let $v_1$ be a factor of $\textbf{v}$ of length $n$, generated by a factor $u_1$ of $\textbf{u}$.

\hspace*{0.5cm} $\bullet$ If $|u_1|=(k+1)q+r+1$, then we have $|v_1|_c=(l-1)q+l-1$ and $v_1$ contains exactly $q+1$ substitutions of $c^l$. Consequently, any substitution of $c^l$ in $u_1$ starts at a position less than or equal to $r$. $v_1$ is written as $v_1=w_0c^lw_1c^lw_2\dotsb c^lw_qc^{l} w_{q+1}$ where $u_1=w_0x_{1}w_1x_{2}w_2\dotsb x_{q}w_qx_{q+1}w_{q+1}$ with $ |w_i|=k$, for $i=1,\ldots,q$ and $0\leq |w_0|,|w_{q+1}|\leq r$. Thus, the factor $u_1$ produces $r+1$ distinct factors of $\textbf{v}$.

\hspace*{0.5cm} $\bullet$ If $|u_1|=(k+1)q+r+j$ with $j=1, \dotsb, k-r-1$ then, we have $|v_1|_c=(l-1)q+l-j$. Consequently, $v_1$ contains $q$ substitutions of $c^l$. Therefore, $u_1$ produces two factors of the form $v_1=w_0c^lw_1c^lw_2\dotsb c^lw_qc^{l-j} $ and $v_2=c^{l-j} t_1c^lt_2\dotsb c^lt_qc^{l}t_{q+1} $ with $u_1=w_0x_1w_1x_2w_2\dotsb x_qw_q =t_1y_1t_2y_2\dotsb t_qy_qt{q+1}$ and $|w_i|=|t_i|=k$, $|w_0|=|t_{q+1}|= r+j$.

\hspace*{0.5cm} $\bullet$ If $|u_1|=(k+1)q+k$, then we have $|v_1|_c=lq+\alpha-k$ and $v_1$ contains $q$ substitutions of $c^l$. Thus, $v_1$ has the form $v_1=c^{h}w_0c^lw_1c^lw_2\dotsb c^lw_qc^{\alpha-k-h}$ with $h=0,1,\ldots, \alpha-k$ and $u_1=w_0x_1w_1x_2w_2\dotsb x_qw_q$ with $|w_i|=k$. Hence, we have:
\begin{align*}
\textbf{P}_\textbf{v}(n)&=(r+1)\textbf{P}_\textbf{u}((k+1)q+r+1)-2\textbf{S}_\textbf{u}((k+1)q+r+1)+2[\textbf{P}_\textbf{u}((k+1)q+r+1)+\dotsb +\\
& \textbf{P}_\textbf{u}((k+1)q+k-1)]+(\alpha-k+1)\textbf{P}_\textbf{u}((k+1)q+k)\\
&=(\alpha-l+1)\textbf{P}_\textbf{u}((k+1)q+\alpha-l+1)-2\textbf{S}_\textbf{u}((k+1)q+\alpha-l)+ 2[\textbf{P}_\textbf{u}((k+1)q+\alpha-l+1)\\
&+\dotsb +
\textbf{P}_\textbf{u}((k+1)q+k-1)]+(\alpha-k+1)\textbf{P}_\textbf{u}((k+1)q+k)\\
&=(\alpha-l+1)\textbf{P}_\textbf{u}((k+1)q+\alpha-l+1)-2\textbf{S}_\textbf{u}((k+1)q+\alpha-l)+\\
& 2\displaystyle\sum_{i=\alpha-l+1}^{k-1} \textbf{P}_\textbf{u}((k+1)q+i)+(\alpha-k+1)\textbf{P}_\textbf{u}((k+1)q+k).
\end{align*}

 \end{itemize} 
\noindent\textbf{Case 2.} Now let us assume $l\leq k$. Then, we also have $m_{0}=l$ and $M_{0}=k$.
\begin{itemize}
\item[(i)] Let us show the following subcases:
\textbf{Case 2.1.} if $n\leq l$. Then, some factors of $\textbf{v}$ of length $n$ contain letter power of $c$. So, they come from the factors of $\textbf{u}$ of length at most $n-$1. These factors of $\textbf{v}$ of length $n$ begin or end with a power of $c$. Moreover, since $n\leq k$ then, any factor of $\textbf{u}$ of length $n$ is also a factor of $\textbf{v}$.
 Thus, we get:
 
 $$L_n(\textbf{v})=L_n(\textbf{u}) \cup\{c^n\hspace{0cm}, \hspace{0.1cm} c^{n-j}t \ tc^{n-j} : \ t\in L_j(\textbf{u}),\ j=1,\ldots, n-1 \}.$$
 
 Hence, we have:
 \begin{align*}
\textbf{P}_\textbf{v}(n)&=\textbf{P}_\textbf{u}(n)+\textbf{P}_u(0)+2[\textbf{P}_\textbf{u}(1)+\textbf{P}_\textbf{u}(2)+\dotsb +\textbf{P}_\textbf{u}(n-1)]\\
&=\textbf{P}_\textbf{u}(n)+1+2\displaystyle\sum_{i=1}^{n-1} \textbf{P}_\textbf{u}(i).
 \end{align*}
  
\textbf{Case 2.2.}  if $l< n\leq k$.
   Then, any factor of $\textbf{u}$ of length $n$ is also a factor of $\textbf{v}$.
 Furthermore, some factors of $\textbf{v}$ of length $n$ contain at most $l$ occurrences of $c$ becaue $n>l$. Then, these factors come from the fators of $\textbf{u}$ of length $n-j$ with $j=1, \cdots, l$. This gives us: $\vspace{0.2cm}$
 
 $L_n(\textbf{v})=L_n(\textbf{u}) \cup \left\{
  \hspace{0.1cm} c^{j}t \hspace{0cm},\hspace{0.1cm} tc^{j}:\hspace{0.1cm}  t \in L_{n-j}(\textbf{u}), j=1, \dotsb , l-1 \right\} \vspace{0.2cm}$ 
 
$\hspace{1.2cm} \cup \left\{ c^lx_2\dotsb x_{n-l+1} , \hspace{0.1cm} x_1c^{l}x_3\dotsb x_{n-l+1} ,\dotsb , \hspace{0.1cm} x_1\dotsb c^{l}x_{n-l+1}, \hspace{0.1cm} x_1\dotsb x_{n-l}c^{l} : \hspace{0.1cm} \hspace{0.1cm} x_1\dotsb x_{n-l+1} \in L(\textbf{u}) \right\}$.$\vspace{0.2cm}$

Thus, we have:$\vspace{0.2cm}$ 
  \begin{align*}
\textbf{P}_\textbf{v}(n)&=\textbf{P}_\textbf{u}(n)+2[\textbf{P}_\textbf{u}(n-l+1)+\dotsb +\textbf{P}_\textbf{u}(n-1)]+(n-l+1)\textbf{P}_\textbf{u}(n-l+1)-2\textbf{S}_\textbf{u}(n)\\
&=\hspace{0cm}(n-l+1)\textbf{P}_\textbf{u}(n-l+1)-2\textbf{S}_\textbf{u}(n)+\textbf{P}_\textbf{u}(n)+2\displaystyle\sum_{i=n-l+1}^{n-1} \textbf{P}_\textbf{u}(i).
  \end{align*}
   
\textbf{Case 2.3.} if $ k<n<k+l$. Since $k<n$, then the factors of $\textbf{v}$ of length $n$ contain $l-1$, $l-2$, $\dotsb$, $n-k+1$ and $n-k$ occurrences of $c$. Consequently, these factors come respectively from the factors of $\textbf{u}$ of length $n-l+1$, $n-l+2$, $\dotsb$, $k-1$ and $k$. This gives us:$\vspace{0.2cm}$
 
 $L_n(\textbf{v})= \left\{ c^lx_2\dotsb x_{n-l+1} , \hspace{0.1cm} x_1c^{l}x_3\dotsb x_{n-l+1} ,\dotsb , \hspace{0.1cm} x_1\dotsb c^{l}x_{n-l+1}, \hspace{0.1cm} x_1\dotsb x_{n-l}c^{l} : \hspace{0.1cm} \hspace{0.1cm} x_1\dotsb x_{n-l+1} \in L(\textbf{u}) \right\} \vspace{0.2cm}$
 
$\hspace{1.2cm} \cup \left\{c^{j}x_1\dotsb x_{n-j} \hspace{0cm},\hspace{0.1cm} x_1\dotsb x_{n-j}c^{j} : \hspace{0.1cm} x_1\dotsb x_{n-j} \in L_{n-j}(u), \hspace{0.1cm} j=n-k, \dotsb , l-1 \right\} \vspace{0.2cm} $

$\hspace{1.2cm}\cup \left\{ c^{j}x_1\dotsb x_kc^{n-k-j}:   \hspace{0.1cm} x_1\dotsb x_k \in L_k(\textbf{u}),  \hspace{0.1cm} j=1, \dotsb , n-k-1 \right\}$.$\vspace{0.2cm}$

 Hence, we have:
 \begin{align*}
 \textbf{P}_\textbf{v}(n)&=(n-l+1)\textbf{P}_\textbf{u}(n-l+1)-2\textbf{S}_\textbf{u}(n-l)+2[\textbf{P}_\textbf{u}(n-l+1)+\dotsb +\textbf{P}_\textbf{u}(k)]+(n-k-1)\textbf{P}_\textbf{u}(k)\\
&= (n-l+1)\textbf{P}_\textbf{u}(n-l+1)-2\textbf{S}_\textbf{u}(n-l)+2[\textbf{P}_\textbf{u}(n-l+1)+\dotsb +\textbf{P}_\textbf{u}(k-1)]+\\
&(n-k+1)\textbf{P}_\textbf{u}(k)\\
&= (n-l+1)\textbf{P}_\textbf{u}(n-l+1)-2\textbf{S}_\textbf{u}(n-l)+(n-k+1)\textbf{P}_\textbf{u}(k)+2\displaystyle\sum_{i=n-l+1}^{k-1} \textbf{P}_\textbf{u}(i).
\end{align*}

 
 \item[(ii)] $n\geq k+l$. Then, there are two integers $\alpha$ and $q$ such that $n=(k+l)q+\alpha$ with $q\geq 1$ and $0\leq \alpha < k+l$. Then, let's reason according to the values of $\alpha$:\\
 
\textbf{Case 2.4.} if $\alpha=0$. By Lemma \ref{(k+i)q}, the factors of $\textbf{v}$ of length $n$ come from the factors of $\textbf{u}$ of length $(k+1)q$. Therefore, these factors contain $(l-1)q$ occurrences of $c$. Thus, any factor has $q$ or $q+1$ substitutions. As a result, we have the set: \\ 
 
 $L_n(\textbf{v})=\lbrace w_0c^lw_1c^l\cdots c^lw_q; |w_0|+|w_q|=k, 1\leq |w_0|,|w_q|\leq k-1; w_0x_1w_1x_2\cdots x_qw_q\in L_{(k+1)q}(\textbf{u}) \rbrace$
 
$\hspace{1.2cm} \bigcup\lbrace c^jw_1c^l\cdots c^lw_qc^{l-j}; j\in \lbrace 0,1,2,\cdots,l-1,l\rbrace; w_1x_2\cdots x_qw_q\in L_{(k+1)q-1}(\textbf{u}) \rbrace$.
 
  \par This gives us the following complexity: 
 
 $$\textbf{P}_\textbf{v}(n)=(k-1)\textbf{P}_\textbf{u}((k+1)q)+(l+1)\textbf{P}_\textbf{u}((k+1)q-1).$$
\textbf{Case 2.5.} if $1\leq \alpha \leq l$. Let $v_1$ be a factor of $\textbf{v}$ of length $n$, generated by a factor $u_1$ of $\textbf{u}$. By Lemma \ref{(k+i)q}, we have $|u_1|=(k+1)q+j$ with $j=0$, 1, $ \dotsb$, $\alpha$.
 
\hspace*{0.5cm} $\bullet $ If $|u_1|=(k+1)q$, then we have $|v_1|_c=(l-1)q+\alpha$. Thus, $v_1$ contains $q-1$ internal substitutions and one substitution at the beginning and another at the end of $v_1$. Hence, $v_1$ can be written in the form $v_1=c^iw_1c^lw_2\dotsb c^lw_qc^{\alpha+l-i} $ where $i=\alpha, \dotsb, l$ with $u_1=x_1w_1x_2w_2\dotsb x_qw_q $ or $u_2=w_1x_1w_2x_2\dotsb w_q x_{q}$ and $|w_i|=k$. Note also that the factors $v_1$ is generated by $w_1x_2w_2\dotsb x_qw_q=x_{1}^{-1}u_1=u_2x_{q}^{-1}$ which is of length $(k+1)q-1$. Consequently, $u_1$ produces ($l-\alpha+1$) distinct factors of $\textbf{v}$ of length $n$.

\hspace*{0.5cm} $\bullet$ If $|u_1|=(k+1)q+j$ with $j=1, \dotsb, \alpha-1$, then we have $|v_1|_c=(l-1)q+\alpha-j$. Thus, $v_1$ contains $q$ substitutions of $c^l$ and one  external substitution of $c^{\alpha-j}$. We can see that the first substitution of $c^l$ in $u_1$ is either at position $j$ or at position $k$. So, $u_1$ produces two factors of the form $v_1=w_0c^lw_1c^lw_2\dotsb c^lw_qc^{\alpha-j} $ and $v_2=c^{\alpha-j} t_1c^lt_2\dotsb c^lt_qc^{l}t_{q+1} $ with $u_1=w_0x_1w_1x_2w_2\dotsb x_qw_q =t_1y_1t_2y_2\dotsb t_qy_qt{q+1}$ and $|w_i|=|t_i|=k$ for $i$=1, 2, $ \dotsb$, $q$ and $|w_0|=|t_{q+1}|= j$.

\hspace*{0.5cm} $\bullet$ If $|u_1|=(k+1)q+\alpha$, then we have $|v_1|_c=(l-1)q$. Thus, $v_1$ contains $q$ internal substitutions of $c^l$. Consequently, the first substitution of $c^l$ in $u_1$ is between the positions $\alpha$ and $k$. Hence, the factor $v_1$ is of the form $v_1=w_0c^lw_1c^lw_2\dotsb c^lw_q$ and $u_1=w_0x_1w_1x_2w_2\dotsb x_qw_q$ with $|w_0|+|w_q|=k+\alpha$, $\alpha \leq |w_0|$, $|w_q|\leq k$. Therefore, $u_1$ produces ($k-\alpha+1$) distinct factors of $\textbf{v}$ of length $(k+l)q+\alpha$. Hence, we obtain:
\begin{align*}
\textbf{P}_\textbf{v}(n)&=(l-\alpha+1)\textbf{P}_\textbf{u}((k+1)q-1)+2[\textbf{P}_\textbf{u}((k+1)q+1)+\dotsb +\textbf{P}_\textbf{u}((k+1)q+\alpha-1)]+\\
&(k-\alpha+1)\textbf{P}_\textbf{u}((k+1)q+\alpha)\\ &=(l-\alpha+1)\textbf{P}_\textbf{u}((k+1)q)+2\displaystyle\sum_{i=1}^{\alpha-1} \textbf{P}_\textbf{u}((k+1)q+i)+(k-\alpha+1)\textbf{P}_\textbf{u}((k+1)q+\alpha.
\end{align*}
   
 \textbf{Case 2.6.} if $l<\alpha\leq k$. Let $v_1$ be a factor of $\textbf{v}$ of length $n$, generated by a factor $u_1$ of $\textbf{u}$. By Lemma \ref{(k+i)q}, we have $|u_1|=(k+1)q+j$ with $j=\alpha-l+1, \cdots , \alpha$.
 
\hspace*{0.5cm} $\bullet $ If $|u_1|=(k+1)q+\alpha-l+1$. Thus, $v_1$ contains $q-1$ internal substitutions and one substitution at the beginning and one at the end of $v_1$. So, $v_1$ can be written as $v_1=w_0c^lw_1c^lw_2\dotsb c^lw_qc^{l}w_{q+1} $ with $|w_i|=k$, $i=1,\dotsb,q$ and $0\leq|w_0|,|w_{q+1}|\leq\alpha-l$.

\hspace*{0.5cm} $\bullet$ If $|u_1|=(k+1)q+j$ with $j=\alpha-l+1, \dotsb, \alpha-1$, then we have $|v_1|_c=(l-1)q+\alpha-j$. Thus, $v_1$ contains $q$ internal substitutions of $c^l$ and one external substitution of $c^{\alpha-j}$. As a result, the first letter substituted from $c^l$ into $u_1$ is either in position $j+1$ or in position $k+1$. So, $u_1$ produces two factors of the form $v_1=w_0c^lw_1c^lw_2\dotsb c^lw_qc^{\alpha-j} $ and $v_2=c^{\alpha-j} t_1c^lt_2\dotsb c^lt_qc^{l}t_{q+1} $ with $u_1=w_0x_1w_1x_2w_2\dotsb x_qw_q =t_1y_1t_2y_2\dotsb t_qy_qt{q+1}$ and $|w_i|=|t_i|=k$, $|w_0|=|t_{q+1}|= j$.

\hspace*{0.5cm}$\bullet$ If $|u_1|=(k+1)q+\alpha$, then we have $|v_1|_c=(l-1)q+k-\alpha$. Thus, $v_1$ contains $q$ substitutions of $c^l$. As a result, $v_1$ is of the form $v_1=c^{h}w_0c^lw_1c^lw_2\dotsb c^lw_qc^{k-\alpha-h}$ with $h=0,1,2,\cdots, k-\alpha$ and $u_1=w_0x_1w_1x_2w_2\dotsb x_qw_q$ with $|w_i|=k$. 

Therefore, $u_1$ produces ($k-\alpha+1$) distinct factors of $\textbf{v}$ of length $(k+l)q+\alpha$. Hence, we obtain:
\begin{align*}
\textbf{P}_\textbf{v}(n)&=(\alpha-l+1)(\textbf{P}_\textbf{u}((k+1)q+\alpha-l+1))-2\textbf{S}_\textbf{u}((k+1)q+\alpha-l)+2[\textbf{P}_\textbf{u}((k+1)q+\alpha-l+1)+\dotsb +\\
&\textbf{P}_\textbf{u}((k+1)q+\alpha-1)]+(k-\alpha+1)\textbf{P}_\textbf{u}((k+1)q+\alpha)\\
&=(\alpha-l+1)(\textbf{P}_\textbf{u}((k+1)q+\alpha-l+1))-2\textbf{S}_\textbf{u}((k+l)q+\alpha-l)+2\displaystyle\sum_{i=\alpha-l+1}^{\alpha-1} \textbf{P}_\textbf{u}((k+1)q+i)+\\
&(k-\alpha+1)\textbf{P}_\textbf{u}((k+1)q+\alpha). 
\end{align*}

\textbf{Case 2.7.} if $\alpha>k$. By Lemma \ref{(k+i)q}, the factors of $\textbf{v}$ of length $n$ come from those of $\textbf{u}$ of length $(k+1)q+\alpha-l+1$, $(k+1)q+\alpha-l+2$, $ \dotsb$, $(k+1)q+k-1$ and $(k+1)q+k$. Let $v_1$ be a factor of $\textbf{v}$ of length $n$, generated by a factor $u_1$ of $\textbf{u}$.

\hspace*{0.5cm} $\bullet$ If $|u_1|=(k+1)q+\alpha-l+1$, then, we have $|v_1|_c=(l-1)q+l-1$ and $v_1$ contains exactly $q+1$ substitutions of $c^l$. Consequently, any substitution of $c^l$ in $u_1$ starts at a position less than or equal to $\alpha-l$. The factor $v_1$ is written as $v_1=w_0c^lw_1c^lw_2\dotsb c^lw_qc^{l} w_{q+1}$ with $0\leq |w_0|,|w_{q+1}|\leq \alpha-l$. Thus, $u_1$ produces $\alpha-l+1$ distinct factors of $\textbf{v}$.

\hspace*{0.5cm} $\bullet$ If $|u_1|=(k+1)q+\alpha-l+j$ with $j=1, \dotsb, k+l-\alpha-1$, then we have $|v_1|_c=(l-1)q+l-j$. As a result, $v_1$ contains $q$ internal substitutions and one substitution at the beginning or one substitution at the end. Therefore, $u_1$ produces two factors of the form $v_1=w_0c^lw_1c^lw_2\dotsb c^lw_qc^{l-j} $ and $v_2=c^{l- j} t_1c^lt_2\dotsb c^lt_qc^{l}t_{q+1} $ with $u_1=w_0x_1w_1x_2w_2\dotsb x_qw_q =t_1y_1t_2y_2\dotsb t_qy_qt{q+1}$ and $|w_i|=|t_i|=k$, $|w_0|=|t_{q+1}|= \alpha-l+j$. Thus any factor $u_1$ produces two distinct factors of $\textbf{v}$.

\hspace*{0.5cm} $\bullet$ If $|u_1|=(k+1)q+k$, then we have $|v_1|_c=lq+\alpha-k$ and $v_1$ contains $q$ internal substitutions of $c^l$. Thus, $v_1$ is of the form $v_1=c^{h}w_0c^lw_1c^lw_2\dotsb c^lw_qc^{\alpha-k-h}$ with $h=0,1,2,\cdots, \alpha-k$ and $|w_i|=k$. Thus, $u_1$ produces $\alpha-k+1$ distinct factors of $\textbf{v}$.

 Hence, we obtain: 
 \begin{align*}
\textbf{P}_\textbf{v}(n)&=(\alpha-l+1)\textbf{P}_\textbf{u}((k+1)q+\alpha-l+1)-2\textbf{S}_\textbf{u}((k+1)q+\alpha-l+1)+\\
&2[\textbf{P}_\textbf{u}((k+1)q+\alpha-l+1)+\dotsb +\textbf{P}_\textbf{u}((k+1)q+k-1)]+
(\alpha-k+1)\textbf{P}_\textbf{u}((k+1)q+k)\\
&=(\alpha-l+1)\textbf{P}_\textbf{u}((k+1)q+\alpha-l+1)-2\textbf{S}_\textbf{u}((k+1)q+\alpha-l+1)+\\
&2[\textbf{P}_\textbf{u}((k+1)q+\alpha-l+1)+\dotsb +\textbf{P}_\textbf{u}((k+1)q+k-1)]+(\alpha-k+1)\textbf{P}_\textbf{u}((k+1)q+k)\\
&=(\alpha-l+1)\textbf{P}_\textbf{u}((k+1)q+\alpha-l+1)-2\textbf{S}_\textbf{u}((k+1)q+\alpha-l+1)+\\
&2\displaystyle\sum_{i=\alpha-l+1}^{k-1} \textbf{P}_\textbf{u}((k+1)q+i)+(\alpha-k+1)\textbf{P}_\textbf{u}((k+1)q+k).
 \end{align*}
 \hfill $\square$
\end{itemize}
 \begin{Corollaire}\label{mod} The substitution does not preserve the modulo-recurrence.
\end{Corollaire}  
\textbf{Proof.} By Theorem \ref{prop-mod}, we have $\textbf{P}_\textbf{v}(k+l)=(k-1)\textbf{P}_\textbf{\emph{u}}(k+1)+(l+1)\textbf{P}_\textbf{\emph{u}}(k) $.

On the other hand, we have $\textbf{P}^f_\textbf{v}(k+l)= \textbf{P}_\textbf{u}(k+1)$. Hence, we obtain $\textbf{P}^f_\textbf{v}(k+l)\neq \textbf{P}_\textbf{v}(k+l)$. By Theorem \ref{theo-stur-mod}, we deduce $\textbf{v}$ is not a modulo-recurrent word.\hfill $\square$
 
 \subsection{ $k$ to $k$ substitution of a letter power in Sturmian words}
In this subsection, we consider a Sturmian word $\textbf{u}$ and we establish the complexity function of the $k$ to $k$ substitution $\textbf{v}=S_{k}^{c^l}(\textbf{u})$. Let us put  $m_0=\min  \left\{l, k\right\}$ and $M_0=\max  \left\{l, k\right\}$.
\begin{theorem}\label{cci}
 
 The complexity function of the word $\textbf{\emph{v}}$ is given  for all integer $n>0$, by:
  $$\hspace{0 cm} \textbf{P}_\textbf{\emph{v}}(n) = \left \{
 \begin{array}{l}
 \hspace{0cm}n^2+2n \hspace*{1.7cm} \text{if} \ n\leq m_0 \\
 \hspace{0cm}(k+1)n+k \hspace{0.9cm} \text{if} \  m_0 < n \leq M_0\ \text{and} \ k < l\\
 \hspace{0cm}n^2+2n-1 \hspace*{1cm} \text{if} \ m_0<n\leq M_0 \ \text{and} \ l<k\\
 \hspace{0cm} (k+1)n+k-1 \ \text{if } \ n> M_0.
 \end{array}
 \right. $$

 \end{theorem}
 
\noindent\textbf{Proof.} By Theorem \ref{stur-mod} $\textbf{u}$ is modulo-recurrent since it is Sturmian. Moreover, $\textbf{u}$ being Sturmian then, we have $\textbf{P}_\textbf{u}(n)=n+1$ and $\textbf{S}_\textbf{u}(n)=1$, for all $n$ according to Definition \ref{def-stur}.
 By applying the Theorem \ref{prop-mod} and reasoning according to the values of $n$, $k$ and $l$, we obtain the different complexities.

\noindent\textbf{Case 1.} If $k \leq l$. Then, we have $m_{0}=k$ and $M_{0}=l$. The following cases arise:\\
\hspace*{0.5cm}\textbf{Case 1.1.} If $n\leq k$. By applying the Theorem \ref{prop-mod}, we obtain:
 \begin{align*}
 \textbf{P}_\textbf{v}(n)&=\textbf{P}_\textbf{u}(n)+\textbf{P}_\textbf{u}(0)+2[\textbf{P}_\textbf{u}(1)+\textbf{P}_\textbf{u}(2)+\dotsb +\textbf{P}_\textbf{u}(n-1)]\\
 &=n+1+1+2(2+3+\dotsb+n) \\ 
 &=n+2+2\dfrac{(n-1)(n+2)}{2}\\
 &=n^2+2n.
 \end{align*}
\hspace*{0.5cm}\textbf{Case 1.2.} If $k<n \leq l$. Then, by applying the Theorem \ref{prop-mod}, we obtain:
  \begin{align*}
  \textbf{P}_\textbf{v}(n)&=\textbf{P}_\textbf{u}(0)+2[\textbf{P}_\textbf{u}(1)+\textbf{P}_\textbf{u}(2)+\dotsb +\textbf{P}_\textbf{u}(k-1)]+(n-k+1)\textbf{P}_\textbf{u}(k)\\
 & = 1+2(2+3+\dotsb+k)+(n-k+1)(k+1)\\ 
 & =1+k^2+3k+kn+n-k^2-2k-1\\
 &=kn+k+n\\
 &=(k+1)n+k.
  \end{align*}
\hspace*{0.5cm}\textbf{Case 1.3.} If $l<n< k+l$. Then, by applying the Theorem \ref{prop-mod}, we obtain:
 \begin{align*}
 \textbf{P}_\textbf{v}(n)&=(n-l+1)\textbf{P}_\textbf{u}(n-l+1)-2+2[\textbf{P}_\textbf{u}(n-l+1)+\dotsb +\textbf{P}_\textbf{u}(k-1)]+(n-k+1)\textbf{P}_\textbf{u}(k)\\
& = (n-l+1)(n-l+1+1)-2+2(n-l+2+n-l+3+\dotsb+k)+(n-k+1)(k+1)\\
 &=kn+k+l+n-l+1-2\\
 &=kn+k+n-1\\
 &=(k+1)n+k-1.
 \end{align*}
\hspace*{0.5cm}\textbf{Case 1.4.} If $n\geq k+l$. Then, there are two integers $\alpha$ and $q$ such that $n=(k+l)q+\alpha$ with $q\geq 1$ and $0\leq \alpha < k+l$. Let's reason according to the values of $\alpha$. \\
\hspace*{1cm}\textbf{Case 1.4.1.} For $\alpha=0 $,  by applying the Theorem \ref{prop-mod}, we obtain the complexity:
 \begin{align*}
 \textbf{P}_\textbf{v}(n)&=(k-1)\textbf{P}_\textbf{u}((k+1)q)+(l+1)\textbf{P}_\textbf{u}((k+1)q-1)\\
 &=(k-1)((k+1)q+1)+(l+1)((k+1)q)\\
 &=(k+l)(k+1)q+k-1\\
 &=(k+1)n+k-1.
 \end{align*}
\hspace*{1cm}\textbf{Case 1.4.2.} For $1\leq \alpha \leq k$. Then, by applying the Theorem \ref{prop-mod}, we get:
  \begin{align*}
  \textbf{P}_\textbf{v}(n)&=(l-\alpha+1)\textbf{P}_\textbf{u}((k+1)q)-1)+2[\textbf{P}_\textbf{u}((k+1)q+1)+\dotsb +\textbf{P}_\textbf{u}((k+1)q+\alpha-1)]+\\
  &(k-\alpha+1)\textbf{P}_\textbf{u}((k+1)q+\alpha)\\
&= (l-\alpha+1)((k+1)q)+2[(k+1)q+2+\dotsb+(k+1)q+\alpha]+(k-\alpha+1)((k+1)q+\alpha+1)\\
&= k((k+l)q+\alpha)+(k+l)q+\alpha+k-1\\
&=kn+n+k-1\\
&=(k+1)n+k-1. 
  \end{align*}
\hspace*{1cm}\textbf{Case 1.4.3.} For $k<\alpha\leq l$. Then, by applying the Theorem \ref{prop-mod}, we get:
 \begin{align*}
\textbf{P}_\textbf{v}(n)&=(l-\alpha+1)(\textbf{P}_\textbf{u}((k+1)q)-1)+2[\textbf{P}_\textbf{u}((k+1)q+1)+\dotsb +\textbf{P}_\textbf{u}((k+1)q+k-1)]+\\
&(\alpha-k+1)\textbf{P}_\textbf{u}((k+1)q+k)\\
&=(l-\alpha+1)((k+1)q)+2[(k+1)q+2+\dotsb+(k+1)q+k]+(\alpha-k+1)(k+1)q+k+1)\\
&=(k+1)((k+l)q+\alpha)+k-1\\
&=(k+1)n+k-1.  
\end{align*}  
\hspace*{1cm}\textbf{Case 1.4.4.} For $\alpha>l$. Then, by applying the Theorem \ref{prop-mod}, we get:
  \begin{align*}
  \textbf{P}_\textbf{v}(n)&=(r+1)\textbf{P}_\textbf{u}((k+1)q+\alpha-l+1)-2+2[\textbf{P}_\textbf{u}((k+1)q+\alpha-l+1)+\dotsb +\textbf{P}_\textbf{u}((k+1)q+k-1)]+\\
  &(\alpha-k+1)\textbf{P}_\textbf{u}((k+1)q+k)\\
&=(k+1)n+k-1.
  \end{align*}
\noindent\textbf{Case 2.} Now assume $l<k$. Let's reason as before according to the values of $k$, $l$ and $n$.\\
\hspace*{0.5cm}\textbf{Case 2.1.} Let $n\leq l$. Then, by applying the Theorem \ref{prop-mod}, we get:
\begin{align*}
\textbf{P}_\textbf{v}(n)&=\textbf{P}_\textbf{u}(n)+\textbf{P}_\textbf{u}(0)+2[\textbf{P}_\textbf{u}(1)+\textbf{P}_\textbf{u}(2)+\dotsb +\textbf{P}_\textbf{u}(n-1)]\\
&=n+1+1+2(2+3+\dotsb+n\\
&=n+2+2\dfrac{(n-1)(n+2)}{2}\\
&=n^2+2n.
\end{align*} 
\hspace*{0.5cm}\textbf{Case 2.2.} Let $l< n\leq k$.  Then, by applying the Theorem \ref{prop-mod}, we get:
 \begin{align*}
\textbf{P}_\textbf{u}(n)&=\textbf{P}_\textbf{u}(n)+2[\textbf{P}_\textbf{u}(n-l+1)+\textbf{P}_\textbf{u}(n-l+2)+\dotsb +\textbf{P}_\textbf{u}(n-1)]+(n-l+1)\textbf{P}_\textbf{u}(n-l+1)-2\\
&= n+1+2(n-l+2+n-l+3+\dotsb+n)+(n-l+1)(n-l+2)-2\\
&=n^2+2n-1.
 \end{align*}
\hspace*{0.5cm}\textbf{Case 2.3.} Let $ k<n<k+l$. Then, by applying the Theorem \ref{prop-mod}, we get:
   \begin{align*}
  \textbf{P}_\textbf{v}(n)&=(n-l+1)\textbf{P}_\textbf{u}(n-l+1)-2+2[\textbf{P}_\textbf{u}(n-l+1)+\dotsb +\textbf{P}_\textbf{u}(k-1)]+(n-k+1)\textbf{P}_\textbf{u}(k)\\
&=(n-l+1)(n-l+1)+2(n-l+2+n-l+3+\dotsb+k)+(n-k+1)(k+1)\\
&=(k+1)n+k-1. 
   \end{align*}
\hspace*{0.5cm}\textbf{Case 2.4.} For $n\geq k+l$. We reason according to the values of $\alpha$. This reasoning is identical to the case where $k\leq l$ and we obtain the same result. \hfill $\square$

\par In $\cite{BBT}$, the following Lemma has been established.
 
 \begin{Lemme}\label{lem1-rap-x}
Let $\textbf{\emph{u}}$ be a non-trivial and uniformly modulo-recurrent word and $\textbf{\emph{v}}=S_{k}^{x}(\textbf{\emph{u}})$. Let $u_1$ and $u_2$ be two distinct factors long enough of $\textbf{\emph{u}}$ and not coming from the same factor by extension. Then, the substitutions of an internal letter $x$ in $u_1$ and $u_2$ give distinct factors of $\textbf{\emph{v}}$.
\end{Lemme}

This result also holds true if we replace $x$ by $x^l$.

\begin{theorem} \label{cc-interne} Let $\textbf{\emph{u}}$ be a Sturmian word over $\left\{a,b\right\}$ and $\textbf{\emph{v}}=S_{k}^{x^l}(\textbf{\emph{u}})$ with $x\in \left\{a,b\right\}$. Then, there exists $n_k$ such that for all $n\geq n_k$, we have:$$\textbf{P}_\textbf{\emph{v}}(n)= (k+1)n+k-1.$$

\end{theorem}

 \textbf{Proof.} According to Theorem \ref{stur-mod}, Sturmian words are non-trivial and uniformly modulo-recurrent words. By Lemma \ref{lem1-rap-x}, there exists an integer $n_k$ such that two factors of $\textbf{u}$ of lengths greater than $n_k$ not originating from the same factor by extension, generate distinct factors of $\textbf{v}$ by $k$ to $k$ substitutions of $x^l$. Let $m_k$ be the minimum length of the factors of $\textbf{v}$ generated by the factors of $\textbf{u}$ of length $n_k$. Then, by applying the  Theorem \ref{cci}, we obtain $\textbf{P}_\textbf{v}(n)=(k+1)n+k-1$, for all $n\geq m_k$.\hfill $\square$

  \end{document}